\documentclass{article}
\usepackage{amssymb}

\input{tcilatex}
\begin{document}

\begin{center}
\textbf{The integral Cauchy problem for generalized Boussinesq equations
with general leading parts}

\textbf{Veli\ B. Shakhmurov}

Department of Mechanical Engineering, Okan University, Akfirat, Tuzla 34959
Istanbul, E-mail: veli.sahmurov@okan.edu.tr;

\textbf{Rishad Shahmurov}

shahmurov@hotmail.com

University of Alabama Tuscaloosa USA, AL 35487

A\textbf{bstract}
\end{center}

In this paper, the integral initial value problems for Boussinesq type
equations are studied. The equation include the general differential
operators. The existence, uniqueness and regularity properties of solution
of these problems are obtained. By choosing differential operators including
in the equation, the regularity properties of the Cauchy problem for
different type of Boussinesg equations are studied.


\begin{center}
\bigskip\ \ \textbf{AMS: 35Lxx, 35Qxx, 47D}

\textbf{1}. \textbf{Introduction}
\end{center}

The aim in this paper is to study the  existence and uniqueness of solution
of the integral initial value problem (IVB) for the general\i zed Boussinesq
equation%
\begin{equation}
u_{tt}+L_{0}u_{tt}+L_{1}u=L_{2}f\left( u\right) ,\text{ }x\in R^{n},\text{ }%
t\in \left( 0,T\right) ,  \tag{1.1}
\end{equation}

\begin{equation}
u\left( x,0\right) =\varphi \left( x\right) +\dint\limits_{0}^{T}\alpha
\left( \sigma \right) u\left( x,\sigma \right) d\sigma ,\text{ }  \tag{1.2}
\end{equation}

\[
u_{t}\left( x,0\right) =\psi \left( x\right) +\dint\limits_{0}^{T}\beta
\left( \sigma \right) u_{t}\left( x,\sigma \right) d\sigma ,
\]%
where $L_{i}$ are d\i fferent\i al operators with constant coefficients, $%
f(u)$ is the given nonlinear function, $\varphi \left( x\right) $ and $\psi
\left( x\right) $ are the given initial value functions, $\alpha $ and $%
\beta $ are measurable functions on $\left( 0,T\right) $.

\textbf{Remark 1.1. }Note that particularly, the condition $\left(
1.2\right) $\ can be expressed as the following multipoint initial condition 
\begin{equation}
u\left( x,0\right) =\varphi \left( x\right) +\dsum\limits_{k=1}^{l}\alpha
_{k}u\left( x,\lambda _{k}\right) ,\text{ }u_{t}\left( 0,x\right) =\psi
\left( x\right) +\dsum\limits_{k=1}^{l}\beta _{k}u_{t}\left( x,\lambda
_{k}\right) .  \tag{1.3}
\end{equation}%
By choosing the operators $L_{i}$ we obtain numerous classes of generalized
Boussinesq type equations which occur in a wide variety of physical systems,
such as in the propagation of longitudinal deformation waves in an elastic
rod, hydro-dynamical process in plasma, in materials science which describe
spinodal decomposition and in the absence of mechanical stresses (see $\left[
1-4\right] $).  For example, if we choose $L_{0}=$ $L_{1}=L_{2}=-\Delta $,
where $\Delta $ is $n-$ dimensioned Laplace, we obtain the Cauchy problem
for the Boussinesq equation%
\begin{equation}
u_{tt}-\Delta u_{tt}-\Delta u=\Delta f\left( u\right) ,\text{ }x\in R^{n},%
\text{ }t\in \left( 0,T\right) ,  \tag{1.4}
\end{equation}%
\[
u\left( x,0\right) =\varphi \left( x\right) ,\text{ }u_{t}\left( x,0\right)
=\psi \left( x\right) .
\]%
The equation $(1.4)$ arises in different situations (see $[1,2]$). For
example, for $n=1$ it describes a limit of a one-dimensional nonlinear
lattice $\left[ 3\right] $, shallow-water waves $\left[ 4,5\right] $ and the
propagation of longitudinal deformation waves in an elastic rod $\left[ 6%
\right] $. Rosenau $\left[ 7\right] $ derived the equations governing
dynamics of one, two and three-dimensional lattices. One of those equations
is $\left( 1.4\right) $. Note that, the existence of solutions and
regularity properties for different type Boussinesq equations are considered
e.g. in $\left[ \text{8-15}\right] .$ In $\left[ 8\right] $ and $\left[ 9%
\right] $ the existence of the global classical solutions and the blow-up of
the solutions of the initial boundary value problem $\left( 1.4\right) $ are
studied. In this paper, we obtain the existence and uniqueness of solution
and regularity properties of the problem $(1.1)-(1.2).$ The strategy is to
express the Boussinesq equation as an integral equation. To treat the
nonlinearity as a small perturbation of the linear part of the equation, the
contraction mapping theorem is used. Also, a priori estimates on $L^{p}$
norms of solutions of the linearized version are utilized. The key step is
the derivation of the uniform estimate of the solutions of the linearized
Boussinesq equation. The methods of harmonic analysis, operator theory,
interpolation of Banach Spaces and embedding theorems in Sobolev spaces are
the main tools implemented to carry out the analysis.

In order to state our results precisely, we introduce some notations and
some function spaces.

\begin{center}
\textbf{Definitions and} \textbf{Background}
\end{center}

Let $E$ be a Banach space. $L^{p}\left( \Omega ;E\right) $ denotes the space
of strongly measurable $E$-valued functions that are defined on the
measurable subset $\Omega \subset R^{n}$ with the norm

\[
\left\Vert f\right\Vert _{L^{p}}=\left\Vert f\right\Vert _{L^{p}\left(
\Omega ;E\right) }=\left( \int\limits_{\Omega }\left\Vert f\left( x\right)
\right\Vert _{E}^{p}dx\right) ^{\frac{1}{p}},\text{ }1\leq p<\infty ,\text{ }
\]

\[
\left\Vert f\right\Vert _{L^{\infty }\left( \Omega \right) }\ =\text{ess}%
\sup\limits_{x\in \Omega }\left\Vert f\left( x\right) \right\Vert _{E}. 
\]

Let $\mathbb{C}$ denote the set of complex numbers.\ For $E=\mathbb{C}$ the $%
L^{p}\left( \Omega ;E\right) $ denotes by $L^{p}\left( \Omega \right) .$

Let $E_{1}$ and $E_{2}$ be two Banach spaces. $\left( E_{1},E_{2}\right)
_{\theta ,p}$ for $\theta \in \left( 0,1\right) ,$ $p\in \left[ 1,\infty %
\right] $ denotes the interpolation spaces defined by $K$-method $\left[ 
\text{17, \S 1.3.2}\right] $.

\ Let $m$ be a positive integer. $W^{m,p}\left( \Omega \right) $ denotes the
Sobolev space, i.e. space of all functions $u\in L^{p}\left( \Omega \right) $
that have the generalized derivatives $\frac{\partial ^{m}u}{\partial
x_{k}^{m}}\in L^{p}\left( \Omega \right) ,$ $1\leq p\leq \infty $ with the
norm 
\[
\ \left\Vert u\right\Vert _{W^{m,p}\left( \Omega \right) }=\left\Vert
u\right\Vert _{L^{p}\left( \Omega \right) }+\sum\limits_{k=1}^{n}\left\Vert 
\frac{\partial ^{m}u}{\partial x_{k}^{m}}\right\Vert _{L^{p}\left( \Omega
\right) }<\infty . 
\]%
\ \ $\ \ $

Let $L^{s,p}\left( R^{n}\right) $, $-\infty <s<\infty $ denotes
Liouville-Sobolev space of order $s$ which is defined as: 
\[
L^{s,p}=L^{s,p}\left( R^{n}\right) =\left( I-\Delta \right) ^{-\frac{s}{2}%
}L^{p}\left( R^{n}\right) 
\]%
with the norm 
\[
\left\Vert u\right\Vert _{L^{s,p}}=\left\Vert \left( I-\Delta \right) ^{%
\frac{s}{2}}u\right\Vert _{L^{p}\left( R^{n}\right) }<\infty . 
\]%
It clear that $L^{0,p}\left( R^{n}\right) =L^{p}\left( R^{n}\right) .$ It is
known that $L^{m,p}\left( R^{n}\right) =W^{m,p}\left( R^{n}\right) $ for the
positive integer $m$ (see e.g. $\left[ \text{18, \S\ 15}\right] $)$.$ Let $%
S\left( R^{n}\right) $ denote Schwartz class, i.e., the space of rapidly
decreasing smooth functions on $R^{n},$ equipped with its usual topology
generated by seminorms. Let $S^{^{\prime }}\left( R^{n}\right) $ denote the
space of all\ continuous linear operators $L:S\left( R^{n}\right)
\rightarrow \mathbb{C},$ equipped with the bounded convergence topology.
Recall $S\left( R^{n}\right) $ is norm dense in $L^{p}\left( R^{n}\right) $
when $1\leq p<\infty .$

Let $1\leq p\leq q<\infty .$ \ A function $\Psi \in L_{\infty }(R^{n})$ is
called a Fourier multiplier from $L_{p}(R^{n})$ to $L_{q}(R^{n})$ if the map 
$B:$ $u\rightarrow F^{-1}\Psi (\xi )Fu$ for $u\in S(R^{n})$ is well defined
and extends to a bounded linear operator

\[
B:L_{p}(R^{n})\rightarrow L_{q}(R^{n}). 
\]
Let $L_{q}^{\ast }\left( E\right) $ denote the space of all $E-$valued
function space such that 
\[
\left\Vert u\right\Vert _{L_{q}^{\ast }\left( E\right) }=\left(
\int\limits_{0}^{\infty }\left\Vert u\left( t\right) \right\Vert _{E}^{q}%
\frac{dt}{t}\right) ^{\frac{1}{q}}<\infty ,\text{ }1\leq q<\infty ,\text{ }%
\left\Vert u\right\Vert _{L_{\infty }^{\ast }\left( E\right) }=\sup_{t\in
\left( 0,\infty \right) }\left\Vert u\left( t\right) \right\Vert _{E}. 
\]

Here, $F$ denote the Fourier transform. Fourier-analytic representation of
Besov spaces on $R^{n}$ is defined as:%
\[
B_{p,q}^{s}\left( R^{n}\right) =\left\{ u\in S^{^{\prime }}\left(
R^{n}\right) ,\right. \text{ } 
\]%
\[
\left\Vert u\right\Vert _{B_{p,q}^{s}\left( R^{n}\right) }=\left\Vert
F^{-1}t^{\varkappa -s}\left( 1+\left\vert \xi \right\vert ^{\frac{\varkappa 
}{2}}\right) e^{-t\left\vert \xi \right\vert ^{2}}Fu\right\Vert
_{L_{q}^{\ast }\left( L^{p}\left( R^{n}\right) \right) }\text{,} 
\]%
\[
\left\vert \xi \right\vert ^{2}=\dsum\limits_{k=1}^{n}\xi _{k}^{2}\text{, }%
\xi =\left( \xi _{1},\xi _{2},...,\xi _{n}\right) ,\left. p\in \left(
1,\infty \right) \text{, }q\in \left[ 1,\infty \right] \text{, }\varkappa
>s\right\} . 
\]

\ It should be note that, the norm of Besov space does not depends on $%
\varkappa $ (see e.g. ( $\left[ \text{17, \S\ 2.3}\right] $). For $p=q$ the
space $B_{p,q}^{s}\left( R^{n}\right) $ will be denoted by $B_{p}^{s}\left(
R^{n}\right) .$

Note that integral conditions for hyperbolic equations were studied e.g. in $%
\left[ \text{21, 22}\right] .$  In a similar way as $\left[ \text{21}\right] 
$ we obtain

\textbf{Lemma 1.1. }Let $C\left( \xi ,t\right) $ be continuous uniformly
bounded in $\xi \in R^{n}$ such that $\left\vert C\left( t\right)
\right\vert \leq 1$ and \ 
\begin{equation}
\left\vert 1+\dint\limits_{0}^{T}\alpha \left( \sigma \right) \beta \left(
\sigma \right) d\sigma \right\vert >\dint\limits_{0}^{T}\left( \left\vert
\alpha \left( \sigma \right) \right\vert +\left\vert \beta \left( \sigma
\right) \right\vert \right) d\sigma .  \tag{1.1}
\end{equation}%
Then the function $O=O\left( \xi \right) $ defined by%
\[
O\left( \xi \right) =1+\dint\limits_{0}^{T}\dint\limits_{0}^{T}C\left( \xi
,\sigma -\tau \right) \alpha \left( \sigma \right) \beta \left( \tau \right)
d\sigma d\tau -\dint\limits_{0}^{T}\left[ \alpha \left( s\right) +\beta
\left( s\right) \right] C\left( \xi ,s\right) ds
\]%
has an uniformly bounded inverse.

Sometimes we use one and the same symbol $C$ without distinction in order to
denote positive constants which may differ from each other even in a single
context. When we want to specify the dependence of such a constant on a
parameter, say $\alpha $, we write $C_{\alpha }$.

The paper is organized as follows: In Section 1, some definitions and
background are given. In Section 2, we obtain the existence of unique
solution and a priory estimates for solution of the linearized problem $%
(1.1)-\left( 1.2\right) .$ In Section 3, we show the existence and
uniqueness of local strong solution of the problem $(1.1)-\left( 1.2\right) $%
. In the Section 4 we show same applications of the problem $(1.1)-\left(
1.2\right) .$

Sometimes we use one and the same symbol $C$ without distinction in order to
denote positive constants which may differ from each other even in a single
context. When we want to specify the dependence of such a constant on a
parameter, say $h$, we write $C_{h}$.

\begin{center}
\textbf{2. Estimates for linearized equation}
\end{center}

In this section, we make the necessary estimates for solutions of integral
IVB for the linearized Boussinesq equation%
\begin{equation}
u_{tt}+L_{0}u_{tt}+L_{1}u=L_{2}g\left( x,t\right) ,\text{ }x\in R^{n},\text{ 
}t\in \left( 0,T\right) ,  \tag{2.1}
\end{equation}%
\begin{equation}
u\left( x,0\right) =\varphi \left( x\right) +\dint\limits_{0}^{T}\alpha
\left( \sigma \right) u\left( x,\sigma \right) d\sigma ,\text{ }  \tag{2.2}
\end{equation}

\[
u_{t}\left( x,0\right) =\psi \left( x\right) +\dint\limits_{0}^{T}\beta
\left( \sigma \right) u_{t}\left( x,\sigma \right) d\sigma , 
\]%
where 
\[
L_{i}u=\dsum\limits_{\left\vert \alpha \right\vert \leq m_{i}}a_{i\alpha
}D^{\alpha }u,\text{ }a_{i\alpha }\in \mathbb{C}\text{, }i=0,1,2, 
\]%
$\alpha =\left( \alpha _{1},\alpha _{2},...,\alpha _{n}\right) ,$ $\alpha
_{k}$ are natural numbers, $\left\vert \alpha \right\vert =$ $%
\dsum\limits_{k=1}^{n}$ $\alpha _{k}$ and $m_{i}$ are positive integers. Let

\[
\text{ }L_{i}\left( \xi \right) =\dsum\limits_{\left\vert \alpha \right\vert
\leq m_{i}}a_{i\alpha }\xi _{1}^{\alpha _{1}}\xi _{2}^{\alpha _{2}}...\xi
_{n}^{\alpha _{n}}\text{, }i=0,1,2. 
\]

Here, 
\[
X_{p}=L^{p}\left( R^{n}\right) \text{, }1\leq p\leq \infty ,\text{ }%
Y^{s,p}=L^{s,p}\left( R^{n}\right) ,\text{ }Y_{1}^{s,p}= 
\]

\[
L^{s,p}\left( R^{n}\right) \cap L^{1}\left( R^{n}\right) \text{, }Y_{\infty
}^{s,p}=L^{s,p}\left( R^{n}\right) \cap L^{\infty }\left( R^{n}\right) ,
\]%
\begin{equation}
Q=Q\left( \xi \right) =L_{1}\left( \xi \right) \left[ 1+L_{0}\left( \xi
\right) \right] ^{-1},\text{ }L\left( \xi \right) =L_{2}\left( \xi \right) %
\left[ 1+L_{0}\left( \xi \right) \right] ^{-1}.  \tag{2.3}
\end{equation}%
\textbf{Condition 2.1. }Let  $\left( 1.1\right) $ holds and $s>\frac{n}{p}$
for $1<p<\infty $. Assume that $L_{1}\left( \xi \right) \neq 0,$ $%
L_{0}\left( \xi \right) \neq -1$\ and there exist a positive constants $M_{1}
$ and $M_{2}$ depend only on $a_{i\alpha }$\ such that 
\[
\left\vert Q^{\frac{1}{2}}\left( \xi \right) \right\vert \leq M_{1}\left(
1+\left\vert \xi \right\vert \right) ^{s-\frac{n}{p}},\text{ }\left\vert
L\left( \xi \right) Q^{-\frac{1}{2}}\left( \xi \right) \right\vert \leq
M_{2}\left( 1+\left\vert \xi \right\vert \right) ^{s-\frac{n}{p}}
\]%
for all $\xi \in R^{n}.$

\textbf{Remark 2.1. }The Condition 2.1 means that there exists positive
constants $M_{1}$ and $M_{2}$ depend only on $a_{\imath \alpha }$ such that 
\[
\left\vert Q^{-\frac{1}{2}}\left( \xi \right) \right\vert =\left[ \left\vert
L_{1}^{-1}\left( \xi \right) +L_{0}\left( \xi \right) L_{1}^{-1}\left( \xi
\right) \right\vert \right] ^{\frac{1}{2}}\leq M_{1}\left( 1+\left\vert \xi
\right\vert \right) ^{\left( s-\frac{n}{p}\right) },\text{ } 
\]

\[
\text{ }\left\vert L_{2}\left( \xi \right) L_{1}^{-1}\left( \xi \right)
\right\vert \left\vert 1+L_{0}\left( \xi \right) \right\vert ^{-\frac{1}{2}%
}\leq M_{2}\left( 1+\left\vert \xi \right\vert \right) ^{\left( s-\frac{n}{p}%
\right) }
\]%
for all $\xi \in R^{n}.$ By Condition 2.1, $L_{1}^{-1}\left( \xi \right) $
and $\left[ 1+L_{0}\left( \xi \right) \right] ^{-\frac{1}{2}}$\ are
uniformly bounded. Therefore, the inequalities $\left( 2.3\right) $ are
satisfied if 
\[
m_{0}-m_{1}\leq 2\left( s-\frac{n}{p}\right) \text{, }m_{2}-m_{1}\leq
2\left( s-\frac{n}{p}\right) 
\]%
respectively, i.e.$\left( 2.3\right) $ is hold trivially if $%
m_{0}=m_{1}=m_{2}.$

First we need the following lemmas

\textbf{Lemma 2.1. }Let the Conditions 2.1 be satisfied.Then problem $\left(
2.1\right) -\left( 2.2\right) $ has a generalized solution.

\textbf{Proof. }By using of Fourier transform we get from $(2.1)-\left(
2.2\right) $:%
\begin{equation}
\hat{u}_{tt}\left( \xi ,t\right) +Q\left( \xi \right) \hat{u}\left( \xi
,t\right) =L\left( \xi \right) \hat{g}\left( \xi ,t\right) ,\text{ } 
\tag{2.4}
\end{equation}%
\begin{equation}
\hat{u}\left( \xi ,0\right) =\hat{\varphi}\left( \xi \right)
+\dint\limits_{0}^{T}\alpha \left( \sigma \right) \hat{u}\left( \xi ,\sigma
\right) d\sigma ,\text{ }  \tag{2.5}
\end{equation}%
\[
\hat{u}_{t}\left( \xi ,0\right) =\hat{\psi}\left( \xi \right)
+\dint\limits_{0}^{T}\beta \left( \sigma \right) \hat{u}_{t}\left( \xi
,\sigma \right) d\sigma , 
\]%
where $\hat{u}\left( \xi ,t\right) $ is a Fourier transform of $u\left(
x,t\right) $ with respect to $x$ and $\hat{\varphi}\left( \xi \right) ,$ $%
\hat{\psi}\left( \xi \right) $ are Fourier transform of $\varphi ,$ $\psi ,$
respectively.

Consider the problem%
\begin{equation}
\hat{u}_{tt}\left( \xi ,t\right) +Q\left( \xi \right) \hat{u}\left( \xi
,t\right) =L\left( \xi \right) \hat{g}\left( \xi ,t\right) ,\text{ } 
\tag{2.6}
\end{equation}%
\[
\hat{u}\left( \xi ,0\right) =u_{0}\left( \xi \right) ,\text{ }\hat{u}%
_{t}\left( \xi ,0\right) =u_{1}\left( \xi \right) ,\text{ }\xi \in R^{n},%
\text{ }t\in \left[ 0,T\right] .\text{ }
\]%
By using the variation of constants it is not hard to see that problem $(2.6)
$ has a unique solution for $\xi \in R^{n}$ and the solution can be expessed
as%
\begin{equation}
\hat{u}\left( \xi ,t\right) =C\left( \xi ,t\right) u_{0}\left( \xi \right)
+S\left( \xi ,t\right) u_{1}\left( \xi \right) +\dint\limits_{0}^{t}S\left(
\xi ,t-\tau \right) \hat{\Phi}\left( \xi ,\tau \right) d\tau ,\text{ }t\in
\left( 0,T\right) ,  \tag{2.7}
\end{equation}%
where, 
\[
C\left( t\right) =C\left( \xi ,t\right) =\cos \left( Q^{\frac{1}{2}}t\right)
,\text{ }S\left( t\right) =S\left( \xi ,t\right) =Q^{-\frac{1}{2}}\sin
\left( Q^{\frac{1}{2}}t\right) ,
\]

\begin{equation}
\hat{\Phi}\left( \xi ,t\right) =L\left( \xi \right) Q^{-\frac{1}{2}}\left(
\xi \right) \sin \left( Q^{\frac{1}{2}}t\right) \hat{g}\left( \xi ,t\right) .
\tag{2.8}
\end{equation}%
\ By using $\left( 2.7\right) $ and the condition 
\[
u_{0}\left( \xi \right) =\hat{\varphi}\left( \xi \right)
+\dint\limits_{0}^{T}\alpha \left( \sigma \right) \hat{u}\left( \xi ,\sigma
\right) d\sigma ,\text{ } 
\]%
we get 
\[
u_{0}\left( \xi \right) =\hat{\varphi}\left( \xi \right)
+\dint\limits_{0}^{T}\alpha \left( \sigma \right) \left[ C\left( \xi ,\sigma
\right) u_{0}\left( \xi \right) +S\left( \xi ,\sigma \right) u_{1}\left( \xi
\right) \right] d\sigma + 
\]

\[
\dint\limits_{0}^{T}\dint\limits_{0}^{\sigma }S\left( \xi ,\sigma -\tau
\right) \hat{\Phi}\left( \xi ,\tau \right) d\tau d\sigma ,\text{ }\tau \in
\left( 0,T\right) . 
\]%
Then, 
\[
\left[ I-\dint\limits_{0}^{T}\alpha \left( \sigma \right) C\left( \xi
,\sigma \right) d\sigma \right] u_{0}\left( \xi \right) -\left[
\dint\limits_{0}^{T}\alpha \left( \sigma \right) S\left( \xi ,\sigma \right)
d\sigma \right] u_{1}\left( \xi \right) = 
\]%
\begin{equation}
\dint\limits_{0}^{T}\dint\limits_{0}^{\sigma }\alpha \left( \sigma \right)
S\left( \xi ,\sigma -\tau \right) \hat{\Phi}\left( \xi ,\tau \right) d\tau
d\sigma +\hat{\varphi}\left( \xi \right) .  \tag{2.9}
\end{equation}

Differentiating both sides of formula $\left( 2.7\right) $ and in view of $%
\left( 2.8\right) $\ we obtain%
\[
\hat{u}_{t}\left( \xi ,t\right) =-Q\left( \xi \right) S\left( \xi ,t\right)
u_{0}\left( \xi \right) +C\left( \xi ,t\right) u_{1}\left( \xi \right) + 
\]%
\begin{equation}
\dint\limits_{0}^{t}C\left( \xi ,t-\tau \right) \hat{\Phi}\left( \xi ,\tau
\right) d\tau ,\text{ }t\in \left( 0,\infty \right) .  \tag{2.10}
\end{equation}%
Using $\left( 2.10\right) $ and the integral condition%
\[
u_{1}\left( \xi \right) =\hat{\psi}\left( \xi \right)
+\dint\limits_{0}^{T}\beta \left( \sigma \right) \hat{u}_{t}\left( \xi
,\sigma \right) d\sigma 
\]%
we obtain%
\[
u_{1}\left( \xi \right) =\hat{\psi}\left( \xi \right)
+\dint\limits_{0}^{T}\beta \left( \sigma \right) \left[ -Q\left( \xi \right)
S\left( \xi ,\sigma \right) u_{0}\left( \xi \right) +C\left( \xi ,\sigma
\right) u_{1}\left( \xi \right) \right] d\sigma + 
\]

\[
\dint\limits_{0}^{T}\dint\limits_{0}^{\sigma }C\left( \xi ,\sigma -\tau
\right) \hat{\Phi}\left( \xi ,\tau \right) d\tau d\sigma . 
\]%
Thus, 
\[
\dint\limits_{0}^{T}\beta \left( \sigma \right) Q\left( \xi \right) S\left(
\xi ,\sigma \right) d\sigma u_{0}\left( \xi \right) +\left[
I-\dint\limits_{0}^{T}\beta \left( \sigma \right) C\left( \xi ,\sigma
\right) d\sigma \right] u_{1}\left( \xi \right) = 
\]%
\begin{equation}
\dint\limits_{0}^{T}\dint\limits_{0}^{\sigma }\beta \left( \sigma \right)
C\left( \xi ,\sigma -\tau \right) \hat{\Phi}\left( \xi ,\tau \right) d\tau
d\sigma +\hat{\psi}\left( \xi \right) .  \tag{2.11}
\end{equation}%
Now, we consider the system of equations $\left( 2.9\right) $, $\left(
2.11\right) $ in $u_{0}\left( \xi \right) $ and $u_{1}\left( \xi \right) $.
The determinant of this system is 
\[
D\left( \xi \right) =\left\vert 
\begin{array}{cc}
\alpha _{11}\left( \xi \right) & \alpha _{12}\left( \xi \right) \\ 
\alpha _{21}\left( \xi \right) & \alpha _{22}\left( \xi \right)%
\end{array}%
\right\vert , 
\]%
where 
\[
\alpha _{11}\left( \xi \right) =I-\dint\limits_{0}^{T}\alpha \left( \sigma
\right) C\left( \xi ,\sigma \right) d\sigma ,\text{ }\alpha _{12}\left( \xi
\right) =-\dint\limits_{0}^{T}\alpha \left( \sigma \right) S\left( \xi
,\sigma \right) d\sigma , 
\]%
\[
\alpha _{21}\left( \xi \right) =\dint\limits_{0}^{T}\beta \left( \sigma
\right) Q\left( \xi \right) S\left( \xi ,\sigma \right) d\sigma ,\text{ }%
\alpha _{22}\left( \xi \right) =I-\dint\limits_{0}^{T}\beta \left( \sigma
\right) S\left( \xi ,\sigma \right) d\sigma . 
\]%
Then by using the properties%
\[
\left[ C\left( \sigma \right) C\left( \tau \right) +Q\left( \xi \right)
S\left( \sigma \right) S\left( \tau \right) \right] =C\left( \sigma -\tau
\right) 
\]
we obtain 
\[
D\left( \xi \right) =I-\dint\limits_{0}^{T}\left[ \alpha \left( \sigma
\right) +\beta \left( \sigma \right) \right] C\left( \sigma \right) d\sigma
+ 
\]%
\[
\dint\limits_{0}^{T}\dint\limits_{0}^{T}\alpha \left( \sigma \right) \beta
\left( \tau \right) \left[ C\left( \xi ,\sigma \right) C\left( \xi ,\tau
\right) +Q\left( \xi \right) S\left( \xi ,\sigma \right) S\left( \xi ,\tau
\right) \right] d\sigma d\tau = 
\]

\[
I-\dint\limits_{0}^{T}\left[ \alpha \left( \sigma \right) +\beta \left(
\sigma \right) \right] C\left( \xi ,\sigma \right) d\sigma
+\dint\limits_{0}^{T}\dint\limits_{0}^{T}C\left( \xi ,\sigma -\tau \right)
\alpha \left( \sigma \right) \beta \left( \tau \right) d\sigma d\tau
=O\left( \xi \right) .
\]%
By Lemma 1.1, $D^{-1}\left( \xi \right) =O^{-1}$ is uniformly bounded.
Solving the system $\left( 2.10\right) -\left( 2.11\right) $, we get%
\begin{equation}
u_{0}\left( \xi \right) =D^{-1}\left( \xi \right) \left\{ \left[
I-\dint\limits_{0}^{T}\beta \left( \sigma \right) C\left( \xi ,\sigma
\right) d\sigma \right] f_{1}\right. +\left. \dint\limits_{0}^{T}\alpha
\left( \sigma \right) S\left( \xi ,\sigma \right) d\sigma f_{2}\right\} , 
\tag{2.12}
\end{equation}%
\[
u_{1}\left( \xi \right) =D^{-1}\left( \xi \right) \left\{ \left[
I-\dint\limits_{0}^{T}\alpha \left( \sigma \right) C\left( \xi ,\sigma
\right) d\sigma \right] f_{2}\right. -\left. \dint\limits_{0}^{T}\left[
\beta \left( \sigma \right) Q\left( \xi \right) S\left( \xi ,\sigma \right)
d\sigma \right] f_{1}\right\} ,
\]%
where 
\[
f_{1}=\dint\limits_{0}^{T}\dint\limits_{0}^{\sigma }\alpha \left( \sigma
\right) S\left( \xi ,\sigma -\tau \right) \hat{\Phi}\left( \xi ,\tau \right)
d\tau d\sigma +\hat{\varphi}\left( \xi \right) ,
\]%
\begin{equation}
f_{2}=\dint\limits_{0}^{T}\dint\limits_{0}^{\sigma }\beta \left( \sigma
\right) C\left( \xi ,\sigma -\tau \right) \hat{\Phi}\left( \tau ,\xi \right)
d\tau d\sigma +\hat{\psi}\left( \xi \right) .  \tag{2.13}
\end{equation}%
From $\left( 2.7\right) ,$ $\left( 2.12\right) $ and $\left( 2.13\right) $
we get that the solution of $\left( 2.4\right) -\left( 2.5\right) $ can be
expressed as 
\[
\hat{u}\left( \xi ,t\right) =D^{-1}\left( \xi \right) \left\{ C\left( \xi
,t\right) \left[ \left( I-\dint\limits_{0}^{T}\beta \left( \sigma \right)
C\left( \xi ,\sigma \right) d\sigma \right) f_{1}\right. \right. +
\]%
\[
\left. \dint\limits_{0}^{T}\alpha \left( \sigma \right) S\left( \xi ,\sigma
\right) d\sigma f_{2}\right] +S\left( t,\xi \right) \left[ \left(
I-\dint\limits_{0}^{T}\alpha \left( \sigma \right) C\left( \xi ,\sigma
\right) d\sigma \right) f_{2}\right. -
\]%
\begin{equation}
\left. \left. \dint\limits_{0}^{T}\beta \left( \sigma \right) Q\left( \xi
\right) S\left( \xi ,\sigma \right) d\sigma f_{1}\right] \right\}
+\dint\limits_{0}^{t}S\left( t-\tau ,\xi \right) \hat{\Phi}\left( \tau ,\xi
\right) d\tau ,\text{ }t\in \left( 0,T\right) .  \tag{2.14}
\end{equation}%
From $\left( 2.14\right) $ we get that there is a generalized solution of $%
(2.1)-(2.2)$ given by 
\begin{equation}
u\left( x,t\right) =S_{1}\left( t\right) \varphi \left( x\right)
+S_{2}\left( t\right) \psi \left( x\right) +\Phi \left( x,t\right) , 
\tag{2.15}
\end{equation}%
where $S_{1}\left( t\right) $ and $S_{2}\left( t\right) $ are defined by 
\[
S_{1}\left( t\right) \varphi =\left( 2\pi \right) ^{-\frac{1}{n}%
}\dint\limits_{R^{n}}\left\{ e^{ix\xi }D^{-1}\left( \xi \right) \right. 
\text{ }
\]%
\[
\left[ C\left( t,\xi \right) \left( I-\dint\limits_{0}^{T}\beta \left(
\sigma \right) C\left( \xi ,\sigma \right) \right)
-\dint\limits_{0}^{T}\beta \left( \sigma \right) L\left( \xi \right) S\left(
\xi ,\sigma \right) \right] d\sigma \left. \hat{\varphi}\left( \xi \right)
d\xi \right\} ,
\]

\begin{equation}
S_{2}\left( t\right) \psi =\left( 2\pi \right) ^{-\frac{1}{n}%
}\dint\limits_{R^{n}}\left\{ e^{ix\xi }D^{-1}\left( \xi \right) C\left(
t,\xi \right) \right.  \tag{2.16}
\end{equation}%
\[
\dint\limits_{0}^{T}\left[ \alpha \left( \sigma \right) S\left( \xi ,\sigma
\right) +S\left( \xi ,\sigma \right) \text{ }\left(
I-\dint\limits_{0}^{T}\alpha \left( \sigma \right) C\left( \xi ,\sigma
\right) \right) d\sigma \right] \left. \hat{\psi}\left( \xi \right) \right\}
d\xi , 
\]

\[
\Phi \left( x,t\right) =\left( 2\pi \right) ^{-\frac{1}{n}%
}\dint\limits_{R^{n}}D^{-1}\left( \xi \right) e^{ix\xi }\left\{
\dint\limits_{0}^{t}S\left( \xi ,t-\tau \right) \Phi \left( \xi ,\tau
\right) d\tau \right. + 
\]

\[
\left[ C\left( \xi ,t\right) \left( I-\dint\limits_{0}^{T}\beta \left(
\sigma \right) C\left( \xi ,\sigma \right) d\sigma \right) \right. + 
\]

\[
\left. S\left( \xi ,t\right) \dint\limits_{0}^{T}\beta \left( \sigma \right)
L\left( \xi \right) S\left( \xi ,\sigma \right) d\sigma \right] g_{1}\left(
\xi \right) +C\left( \xi ,t\right) \dint\limits_{0}^{T}\alpha \left( \sigma
\right) S\left( \xi ,\sigma \right) d\sigma + 
\]%
\[
S\left( \xi ,t\right) \left( I-\dint\limits_{0}^{T}\alpha \left( \sigma
\right) C\left( \xi ,\sigma \right) d\sigma \right) \left. g_{2}\left( \xi
\right) \right\} d\xi , 
\]%
here 
\begin{equation}
g_{1}\left( \xi \right) =\dint\limits_{0}^{T}\dint\limits_{0}^{\sigma
}\alpha \left( \sigma \right) S\left( \xi ,\sigma -\tau \right) \hat{g}%
\left( \xi ,\tau \right) d\tau d\sigma ,  \tag{2.14}
\end{equation}%
\[
g_{2}\left( \xi \right) =\dint\limits_{0}^{T}\dint\limits_{0}^{\sigma }\beta
\left( \sigma \right) C\left( \xi ,\sigma -\tau \right) \hat{g}\left( \xi
,\tau \right) d\tau d\sigma . 
\]

\bigskip \textbf{Theorem 2.1. }Let the Condition 2.1\ be hold. Then for $%
\varphi ,$ $\psi ,$ $g\left( x,t\right) \in Y_{1}^{s,p}$ the solution $%
\left( 2.1\right) -\left( 2.2\right) $ satisfies the following estimate 
\begin{equation}
\left\Vert u\right\Vert _{X_{\infty }}+\left\Vert u_{t}\right\Vert
_{X_{\infty }}\leq C\left[ \left\Vert \varphi \right\Vert
_{Y^{s,p}}+\left\Vert \varphi \right\Vert _{X_{1}}\right. +  \tag{2.15}
\end{equation}

\[
\left\Vert \psi \right\Vert _{Y^{s,p}}+\left\Vert \psi \right\Vert
_{X_{1}}+\left. \dint\limits_{0}^{t}\left( \left\Vert g\left( .,\tau \right)
\right\Vert _{Y^{s,p}}+\left\Vert g\left( .,\tau \right) \right\Vert
_{X_{1}}\right) d\tau \right] 
\]%
uniformly with respect to $t\in \left[ 0,T\right] .$

\textbf{Proof. }Let $N\in \mathbb{N}$ and 
\[
\Pi _{N}=\left\{ \xi :\xi \in R^{n},\text{ }\left\vert \xi \right\vert \leq
N\right\} ,\text{ }\Pi _{N}^{\prime }=\left\{ \xi :\xi \in R^{n},\text{ }%
\left\vert \xi \right\vert \geq N\right\} . 
\]%
It is clear to see that

\[
\left\Vert F^{-1}C\left( \xi ,t\right) \hat{\varphi}\left( \xi \right)
\right\Vert _{X_{\infty }}+\left\Vert F^{-1}S\left( \xi \right) \hat{\psi}%
\left( \xi ,t\right) \right\Vert _{X_{\infty }}\leq 
\]%
\begin{equation}
\left\Vert \dint\limits_{R^{n}}e^{ix\xi }C\left( \xi ,t\right) \varphi
\left( x\right) dx\right\Vert _{L^{\infty }\left( \Pi _{N}\right)
}+\left\Vert \dint\limits_{R^{n}}e^{ix\xi }S\left( \xi ,t\right) \psi \left(
x\right) dx\right\Vert _{L^{\infty }\left( \Pi _{N}\right) }+  \tag{2.16}
\end{equation}%
\[
\left\Vert F^{-1}C\left( \xi ,t\right) \hat{\varphi}\left( \xi \right)
\right\Vert _{L^{\infty }\left( \Pi _{N}^{\prime }\right) }+\left\Vert
F^{-1}S\left( \xi ,t\right) \hat{\psi}\left( \xi \right) \right\Vert
_{L^{\infty }\left( \Pi _{N}^{\prime }\right) }. 
\]%
By using the Minkowski's inequality for integrals and in view of the
uniformly boundedness of $C\left( \xi ,t\right) $, $S\left( \xi ,t\right) $
on $\Pi _{N}$ we have 
\begin{equation}
\left\Vert \dint\limits_{R^{n}}e^{ix\xi }C\left( \xi ,t\right) \varphi
\left( x\right) dx\right\Vert _{L^{\infty }\left( \Pi _{N}\right)
}+\left\Vert \dint\limits_{R^{n}}e^{ix\xi }S\left( \xi ,t\right) \psi \left(
x\right) dx\right\Vert _{L^{\infty }\left( \Pi _{N}\right) }\leq  \tag{2.17}
\end{equation}

\[
C\left[ \left\Vert \varphi \right\Vert _{X_{1}}+\left\Vert \psi \right\Vert
_{X_{1}}\right] . 
\]%
Hence, 
\[
\left\Vert F^{-1}C\left( \xi ,t\right) \hat{\varphi}\left( \xi \right)
\right\Vert _{L^{\infty }\left( \Pi _{N}^{\prime }\right) }+\left\Vert
F^{-1}S\left( \xi ,t\right) \hat{\psi}\left( \xi \right) \right\Vert
_{L^{\infty }\left( \Pi _{N}^{\prime }\right) }= 
\]%
\begin{equation}
=\left\Vert F^{-1}\left( 1+\left\vert \xi \right\vert ^{2}\right) ^{-\frac{s%
}{2}}C\left( \xi ,t\right) \left( 1+\left\vert \xi \right\vert ^{2}\right) ^{%
\frac{s}{2}}\hat{\varphi}\left( \xi \right) \right\Vert _{L^{\infty }\left(
\Pi _{N}^{\prime }\right) }+  \tag{2.18}
\end{equation}%
\[
\left\Vert F^{-1}\left( 1+\left\vert \xi \right\vert ^{2}\right)
^{-s}S\left( \xi ,t\right) \left( 1+\left\vert \xi \right\vert \right) ^{%
\frac{s}{2}}\hat{\psi}\left( \xi \right) \right\Vert _{L^{\infty }\left( \Pi
_{N}^{\prime }\right) }. 
\]%
By using $\left( 2.3\right) $ and\ the first estimate in Condition 2.1 we get%
\[
\sup\limits_{\xi \in R^{n},t\in \left[ 0,T\right] }\left\vert \xi
\right\vert \left\vert ^{\left\vert \alpha \right\vert +\frac{n}{p}%
}D^{\alpha }\left[ \left( 1+\left\vert \xi \right\vert ^{2}\right) ^{-\frac{s%
}{2}}C\left( \xi ,t\right) \right] \right\vert \leq C_{2}, 
\]%
\ 
\begin{equation}
\sup\limits_{\xi \in R^{n},t\in \left[ 0,T\right] }\left\vert \xi
\right\vert \left\vert ^{\left\vert \alpha \right\vert +\frac{n}{p}%
}D^{\alpha }\left[ \left( 1+\left\vert \xi \right\vert ^{2}\right) ^{-\frac{s%
}{2}}S\left( \xi ,t\right) \right] \right\vert \leq C_{2},  \tag{2.19}
\end{equation}%
for $s>\frac{n}{p},$ $\alpha =\left( \alpha _{1},\alpha _{2},...,\alpha
_{n}\right) $, $\alpha _{k}\in \left\{ 0,1\right\} $, $\xi \in R^{n}$ and $%
\xi \neq 0$ uniformly in $t\in \left[ 0,T\right] .$ By multiplier theorems
(see e.g. $\left[ 16\right] $) from $\left( 2.19\right) $ we get that the
functions $\left( 1+\left\vert \xi \right\vert ^{2}\right) ^{-\frac{s}{2}%
}C\left( \xi ,t\right) ,$ $\left( 1+\left\vert \xi \right\vert ^{2}\right)
^{-\frac{s}{2}}S\left( \xi ,t\right) $ are $L^{p}\left( R^{n}\right)
\rightarrow L^{\infty }\left( R^{n}\right) $ Fourier multipliers. Then by
Minkowski's inequality for integrals, from $\left( 2.17\right) -\left(
2.18\right) $ we obtain%
\begin{equation}
\left\Vert F^{-1}C\left( \xi ,t\right) \hat{\varphi}\left( \xi \right)
\right\Vert _{L^{\infty }\left( \Pi _{N}^{\prime }\right) }+\left\Vert
F^{-1}S\left( \xi ,t\right) \hat{\psi}\left( \xi \right) \right\Vert
_{L^{\infty }\left( \Pi _{N}^{\prime }\right) }\leq  \tag{2.20}
\end{equation}

\[
C\left[ \left\Vert \varphi \right\Vert _{Y^{s,p}}+\left\Vert \psi
\right\Vert _{Y^{s,p}}\right] . 
\]

By using the representation of $\hat{\Phi}\left( \xi ,t\right) $ in $\left(
2.8\right) $ and\ the second inequality in Condition 2.1 we get the uniforum
estimate \ 
\begin{equation}
\sup\limits_{\xi \in R^{n},t\in \left[ 0,T\right] }\left\vert \xi
\right\vert \left\vert ^{\left\vert \alpha \right\vert +\frac{n}{p}%
}D^{\alpha }\left[ \left( 1+\left\vert \xi \right\vert ^{2}\right) ^{-\frac{s%
}{2}}\hat{\Phi}\left( \xi ,t\right) \right] \right\vert \leq C_{3}. 
\tag{2.21}
\end{equation}

By reasoning as the above we have 
\[
\left\Vert F^{-1}\dint\limits_{0}^{t}\hat{\Phi}\left( t-\tau ,\xi \right) 
\hat{g}\left( \xi ,\tau \right) d\tau \right\Vert _{X_{\infty }}\leq
C\dint\limits_{0}^{t}\left( \left\Vert g\left( .,\tau \right) \right\Vert
_{Y^{s}}+\left\Vert g\left( .,\tau \right) \right\Vert _{X_{1}}\right) d\tau
. 
\]

Hence, we obtain the estimate%
\begin{equation}
\left\Vert u\right\Vert _{X_{\infty }}\leq C\left[ \left\Vert \varphi
\right\Vert _{Y^{s,p}}+\left\Vert \varphi \right\Vert _{X_{1}}\right. + 
\tag{2.22}
\end{equation}

\[
\left\Vert \psi \right\Vert _{Y^{s,p}}+\left\Vert \psi \right\Vert
_{X_{1}}+\left. \dint\limits_{0}^{t}\left( \left\Vert g\left( .,\tau \right)
\right\Vert _{Y^{s,p}}+\left\Vert g\left( .,\tau \right) \right\Vert
_{X_{1}}\right) d\tau \right] . 
\]

By using $\left( 2.3\right) $ and Condition 2.1 in view of $\left(
2.20\right) $ in similar way, we deduced the estimate of type $\left(
2.22\right) $ for $u_{t}$, i.e. we obtain the assertion.

\textbf{Theorem 2.2. }Let the Conditions 2.1 be hold. Then for $\varphi ,$ $%
\psi ,$ $g\left( x,t\right) \in Y^{s,p}$ the solution of the problem $\left(
2.1\right) -\left( 2.2\right) $ satisfies the following uniform estimate%
\begin{equation}
\left( \left\Vert u\right\Vert _{Y^{s,p}}+\left\Vert u_{t}\right\Vert
_{Y^{s,p}}\right) \leq C\left( \left\Vert \varphi \right\Vert
_{Y^{s,p}}+\left\Vert \psi \right\Vert
_{Y^{s,p}}+\dint\limits_{0}^{t}\left\Vert g\left( .,\tau \right) \right\Vert
_{Y^{s,p}}d\tau \right) .  \tag{2.23}
\end{equation}

\textbf{Proof. }From $\left( 2.7\right) $ we have the following uniform
estimate 
\begin{equation}
\left( \left\Vert F^{-1}\left( 1+\left\vert \xi \right\vert ^{2}\right) ^{%
\frac{s}{2}}\hat{u}\right\Vert _{X_{p}}+\left\Vert F^{-1}\left( 1+\left\vert
\xi \right\vert ^{2}\right) ^{\frac{s}{2}}\hat{u}_{t}\right\Vert
_{X_{p}}\right) \leq  \tag{2.24}
\end{equation}

\[
C\left\{ \left\Vert F^{-1}\left( 1+\left\vert \xi \right\vert \right) ^{%
\frac{s}{2}}C\left( \xi ,t\right) \hat{\varphi}\right\Vert _{X_{p}}\right.
+\left\Vert F^{-1}\left( 1+\left\vert \xi \right\vert \right) ^{\frac{s}{2}%
}S\left( \xi ,t\right) \hat{\psi}\right\Vert _{X_{p}}+ 
\]

\[
\left. \dint\limits_{0}^{t}\left\Vert \left( 1+\left\vert \xi \right\vert
\right) ^{\frac{s}{2}}\hat{\Phi}\left( \xi ,t-\tau \right) \hat{g}\left(
.,\tau \right) \right\Vert _{X_{p}}d\tau \right\} .
\]

\bigskip By Condition 2.1 and by virtue of Fourier multiplier theorems (see $%
\left[ \text{17, \S\ 2.2}\right] $)\ we get that $C\left( \xi ,t\right) $, $%
S\left( \xi ,t\right) $ and $\hat{\Phi}\left( \xi ,t\right) $\ are Fourier
multipliers in $L^{p}\left( R^{n}\right) $ uniformly with respect to $t\in %
\left[ 0,T\right] .$ So, the estimate $\left( 2.24\right) $ by using the
Minkowski's inequality for integrals implies $\left( 2.23\right) .$

\begin{center}
\textbf{3. Initial value problem for nonlinear equation}
\end{center}

In this section, we will show the local existence and uniqueness of solution
for the Cauchy problem $(1.1)-(1.2).$ For the study of the nonlinear problem 
$\left( 1.1\right) -\left( 1.2\right) $ we need the following lemmas

\textbf{Lemma 3.1} (Nirenberg's inequality) $\left[ 19\right] $. Assume that 
$u\in L^{p}\left( \Omega \right) $, $D^{m}u$ $\in L^{q}\left( \Omega \right) 
$, $p,q\in \left( 1,\infty \right) $. Then for $i$ with $0\leq i\leq m,$ $m>%
\frac{n}{q}$ we have 
\begin{equation}
\left\Vert D^{i}u\right\Vert _{r}\leq C\left\Vert u\right\Vert _{p}^{1-\mu
}\dsum\limits_{k=1}^{n}\left\Vert D_{k}^{m}u\right\Vert _{q}^{\mu }, 
\tag{3.1}
\end{equation}%
where%
\[
\frac{1}{r}=\frac{i}{m}+\mu \left( \frac{1}{q}-\frac{m}{n}\right) +\left(
1-\mu \right) \frac{1}{p},\text{ }\frac{i}{m}\leq \mu \leq 1. 
\]

\textbf{Lemma 3.2 }$\left[ 20\right] .$\textbf{\ }Assume that $u\in $ $%
W^{m,p}\left( \Omega \right) \cap L^{\infty }\left( \Omega \right) $ and $%
f\left( u\right) $ possesses continuous derivatives up to order $m\geq 1$.
Then $f\left( u\right) -f\left( 0\right) \in W^{m,p}\left( \Omega \right) $
and 
\[
\left\Vert f\left( u\right) -f\left( 0\right) \right\Vert _{p}\leq
\left\Vert f^{^{\left( 1\right) }}\left( u\right) \right\Vert _{\infty
}\left\Vert u\right\Vert _{p}, 
\]

\begin{equation}
\left\Vert D^{k}f\left( u\right) \right\Vert _{p}\leq
C_{0}\dsum\limits_{j=1}^{k}\left\Vert f^{\left( j\right) }\left( u\right)
\right\Vert _{\infty }\left\Vert u\right\Vert _{\infty }^{j-1}\left\Vert
D^{k}u\right\Vert _{p}\text{, }1\leq k\leq m,  \tag{3.2}
\end{equation}%
where $C_{0}$ $\geq 1$ is a constant.

Let%
\[
\text{ }X_{p}=L^{p}\left( R^{n}\right) ,\text{ }Y=W^{2,p}\left( R^{n}\right)
,\text{ }E_{0}=\left( X_{p},Y\right) _{\frac{1}{2p},p}=B_{p}^{2-\frac{1}{p}%
}\left( R^{n}\right) . 
\]
\textbf{Remark 3.1. }By using J.Lions-I. Petree result (see e.g $\left[ 
\text{17, \S\ 1.8.}\right] $) we obtain that the map $u\rightarrow u\left(
t_{0}\right) $, $t_{0}\in \left[ 0,T\right] $ is continuous and surjective
from $W^{2,p}\left( 0,T\right) $ onto $E_{0}$ and there is a constant $C_{1}$
such that 
\[
\left\Vert u\left( t_{0}\right) \right\Vert _{E_{0}}\leq C_{1}\left\Vert
u\right\Vert _{W^{2,p}\left( 0,T\right) },\text{ }1\leq p\leq \infty \text{.}
\]

First all of, we define the space $Y\left( T\right) =C\left( \left[ 0,T%
\right] ;Y_{\infty }^{2,p}\right) $ equipped with the norm defined by%
\[
\left\Vert u\right\Vert _{Y\left( T\right) }=\max\limits_{t\in \left[ 0,T%
\right] }\left\Vert u\right\Vert _{Y^{2,p}}+\max\limits_{t\in \left[ 0,T%
\right] }\left\Vert u\right\Vert _{X_{\infty }},\text{ }u\in Y\left(
T\right) . 
\]

It is easy to see that $Y\left( T\right) $ is a Banach space. For $\varphi $%
, $\psi \in Y^{2,p}$, let 
\[
M=\left\Vert \varphi \right\Vert _{Y^{2,p}}+\left\Vert \varphi \right\Vert
_{X_{\infty }}+\left\Vert \psi \right\Vert _{Y^{2,p}}+\left\Vert \psi
\right\Vert _{X_{\infty }}. 
\]

\textbf{Definition 3.1. }For any $T>0$ if $\upsilon ,$ $\psi \in Y_{\infty
}^{2,p}$ and $u$ $\in C\left( \left[ 0,T\right] ;Y_{\infty }^{2,p}\right) $
satisfies the equation $(1.1)-(1.2)$ then $u\left( x,t\right) $ is called
the continuous solution\ or the strong solution of the problem $(1.1)-(1.2).$
If $T<\infty $, then $u\left( x,t\right) $ is called the local strong
solution of the problem $(1.1)-(1.2).$ If $T=\infty $, then $u\left(
x,t\right) $ is called the global strong solution of the problem $%
(1.1)-(1.2) $.

\textbf{Condition 3.1. }Assume:

(1) The Condition 2.1 holds,  $\varphi ,$ $\psi $ $\in Y_{\infty }^{2,p}$
for $1<p<\infty $ $\ $and $\frac{n}{p}<2$;

(2) the function $u\rightarrow $ $f\left( x,t,u\right) $: $R^{n}\times \left[
0,T\right] \times E_{0}\rightarrow E$ is a measurable in $\left( x,t\right)
\in R^{n}\times \left[ 0,T\right] $ for $u\in E_{0};$ $f\left( x,t,u\right) $%
. Moreover, $f\left( x,t,u\right) $ is continuous in $u\in E_{0}$ and $%
f\left( x,t,u\right) \in C^{\left( 3\right) }\left( E_{0};E\right) $
uniformly with respect to $x\in R^{n},$ $t\in \left[ 0,T\right] .$ Main aim
of this section is to prove the following result:

\textbf{Theorem 3.1. }Let the Condition 3.1 hold. Then problem $\left(
1.1\right) -\left( 1.2\right) $ has a unique local strange solution $u\in
C^{\left( 2\right) }\left( \left[ 0,\right. \left. T_{0}\right) ;Y_{\infty
}^{2,p}\right) $, where $T_{0}$ is a maximal time interval that is
appropriately small relative to $M$. Moreover, if

\begin{equation}
\sup_{t\in \left[ 0\right. ,\left. T_{0}\right) }\left( \left\Vert
u\right\Vert _{Y^{2,p}}+\left\Vert u\right\Vert _{X_{\infty }}+\left\Vert
u_{t}\right\Vert _{Y^{2,p}}+\left\Vert u_{t}\right\Vert _{X_{\infty
}}\right) <\infty  \tag{3.3}
\end{equation}%
then $T_{0}=\infty .$

\textbf{Proof. }First, we are going to prove the existence and the
uniqueness of the local continuous solution of the problem $(1.1)-\left(
1.2\right) $ by contraction mapping principle. Consider a map $G$ on $%
Y\left( T\right) $ such that $G(u)$ is the solution of the Cauchy problem%
\begin{equation}
G_{tt}\left( u\right) +L_{0}G_{tt}\left( u\right) +L_{1}G\left( u\right)
=L_{2}f\left( G\left( u\right) \right) ,\text{ }x\in R^{n},\text{ }t\in
\left( 0,T\right) ,  \tag{3.4}
\end{equation}%
\[
G\left( u\right) \left( x,0\right) =\varphi \left( x\right)
+\dint\limits_{0}^{T}\alpha \left( \sigma \right) G\left( u\right) \left(
x,\sigma \right) d\sigma ,\text{ } 
\]

\[
G\left( u\right) _{t}\left( x,0\right) =\psi \left( x\right)
+\dint\limits_{0}^{T}\beta \left( \sigma \right) G\left( u\right) _{t}\left(
x,\sigma \right) d\sigma .
\]%
From Lemma 3.2 we know that $f(u)\in $ $L^{p}\left( 0,T;Y_{\infty
}^{2,p}\right) $ for any $T>0$. Thus, by Theorem 2.1, problem $\left(
3.4\right) $ has a unique solution which can be written as%
\[
G\left( u\right) \left( x,t\right) =S_{1}\left( t\right) \varphi \left(
x\right) +S_{2}\left( t\right) \psi \left( x\right) +
\]%
\begin{equation}
+\dint\limits_{0}^{t}F^{-1}\left[ S\left( t-\tau ,\xi \right) L\left( \xi
\right) \hat{f}\left( u\right) \left( \xi ,\tau \right) \right] d\tau ,\text{
}t\in \left( 0,T\right) ,  \tag{3.5}
\end{equation}%
where $S_{1}\left( t\right) $, $S_{2}\left( t\right) $ are linear operators
in $L^{p}\left( R^{n}\right) $ defined by $\left( 2.15\right) $.\ From Lemma
3.2 it is easy to see that the map $G$ is well defined for $f\in C^{\left(
2\right) }\left( X_{0};\mathbb{C}\right) $. We put 
\[
Q\left( M;T\right) =\left\{ u\mid u\in Y\left( T\right) \text{, }\left\Vert
u\right\Vert _{Y\left( T\right) }\leq M+1\right\} .
\]

First, by reasoning as in $\left[ 9\right] $\ let us prove that the map $G$
has a unique fixed point in $Q\left( M;T\right) .$ For this aim, it is
sufficient to show that the operator $G$ maps $Q\left( M;T\right) $ into $%
Q\left( M;T\right) $ and $G:$ $Q\left( M;T\right) $ $\rightarrow $ $Q\left(
M;T\right) $ is strictly contractive if $T$ is appropriately small relative
to $M.$ Consider the function \ $\bar{f}\left( \xi \right) $: $\left[
0,\right. $ $\left. \infty \right) \rightarrow \left[ 0,\right. $ $\left.
\infty \right) $ defined by 
\[
\ \bar{f}\left( \xi \right) =\max\limits_{\left\vert x\right\vert \leq \xi
}\left\{ \left\Vert f^{\left( 1\right) }\left( x\right) \right\Vert _{%
\mathbb{C}},\left\Vert f^{\left( 2\right) }\left( x\right) \right\Vert _{%
\mathbb{C}}\text{ }\right\} ,\text{ }\xi \geq 0. 
\]

It is clear to see that the function $\bar{f}\left( \xi \right) $ is
continuous and nondecreasing on $\left[ 0,\right. $ $\left. \infty \right) .$
From Lemma 3.2 we have\qquad

\[
\left\Vert f\left( u\right) \right\Vert _{Y^{2,p}}\leq \left\Vert f^{\left(
1\right) }\left( u\right) \right\Vert _{X_{\infty }}\left\Vert u\right\Vert
_{X_{p}}+\left\Vert f^{\left( 1\right) }\left( u\right) \right\Vert
_{X_{\infty }}\left\Vert Du\right\Vert _{X_{p}}+ 
\]

\begin{equation}
C_{0}\left[ \left\Vert f^{\left( 1\right) }\left( u\right) \right\Vert
_{X_{\infty }}\left\Vert u\right\Vert _{X_{p}}+\left\Vert f^{\left( 2\right)
}\left( u\right) \right\Vert _{X_{\infty }}\left\Vert u\right\Vert
_{X_{\infty }}\left\Vert D^{2}u\right\Vert _{X_{p}}\right] \leq  \tag{3.6}
\end{equation}

\ 
\[
2C_{0}\bar{f}\left( M+1\right) \left( M+1\right) \left\Vert u\right\Vert
_{Y^{2,p}}.
\]%
By using the Theorem 2.1 we obtain from $\left( 3.5\right) $:%
\begin{equation}
\left\Vert G\left( u\right) \right\Vert _{X_{\infty }}\leq \left\Vert
\varphi \right\Vert _{X_{\infty }}+\left\Vert \psi \right\Vert _{X_{\infty
}}+\dint\limits_{0}^{t}\left\Vert f\left( x,\tau ,u\left( \tau \right)
\right) \right\Vert _{X_{\infty }},  \tag{3.7}
\end{equation}%
\begin{equation}
\left\Vert G\left( u\right) \right\Vert _{Y^{2,p}}\leq \left\Vert \varphi
\right\Vert _{Y^{2,p}}+\left\Vert \psi \right\Vert
_{Y^{2,p}}+\dint\limits_{0}^{t}\left\Vert f\left( x,\tau ,u\left( \tau
\right) \right) \right\Vert _{Y^{2,p}}d\tau .  \tag{3.8}
\end{equation}%
Thus, from $\left( 3.6\right) -\left( 3.8\right) $ and Lemma 3.2 we get 
\[
\left\Vert G\left( u\right) \right\Vert _{Y\left( T\right) }\leq M+T\left(
M+1\right) \left[ 1+2C_{0}\left( M+1\right) \bar{f}\left( M+1\right) \right]
.
\]%
If $T$ satisfies 
\begin{equation}
T\leq \left\{ \left( M+1\right) \left[ 1+2C_{0}\left( M+1\right) \bar{f}%
\left( M+1\right) \right] \right\} ^{-1},  \tag{3.9}
\end{equation}%
then 
\[
\left\Vert Gu\right\Vert _{Y\left( T\right) }\leq M+1.
\]%
Therefore, if $\left( 3.9\right) $ holds, then $G$ maps $Q\left( M;T\right) $
into $Q\left( M;T\right) .$ Now, we are going to prove that the map $G$ is
strictly contractive. Assume $T>0$ and $u_{1},$ $u_{2}\in $ $Q\left(
M;T\right) $ given. We get%
\[
G\left( u_{1}\right) -G\left( u_{2}\right)
=\dint\limits_{0}^{t}F^{-1}S\left( t-\tau ,\xi \right) L\left( \xi \right) 
\left[ \hat{f}\left( u_{1}\right) \left( \xi ,\tau \right) -\hat{f}\left(
u_{2}\right) \left( \xi ,\tau \right) \right] d\tau ,\text{ }t\in \left(
0,T\right) .
\]

By using the assumption (2) and the mean value theorem, we obtain%
\[
\hat{f}\left( u_{1}\right) -\hat{f}\left( u_{2}\right) =\hat{f}^{\left(
1\right) }\left( u_{2}+\eta _{1}\left( u_{1}-u_{2}\right) \right) \left(
u_{1}-u_{2}\right) ,\text{ }
\]

\[
D_{\xi }\left[ \hat{f}\left( u_{1}\right) -\hat{f}\left( u_{2}\right) \right]
=\hat{f}^{\left( 2\right) }\left( u_{2}+\eta _{2}\left( u_{1}-u_{2}\right)
\right) \left( u_{1}-u_{2}\right) D_{\xi }u_{1}+\text{ }
\]%
\[
\hat{f}^{\left( 1\right) }\left( u_{2}\right) \left( D_{\xi }u_{1}-D_{\xi
}u_{2}\right) ,
\]%
\[
D_{\xi }^{2}\left[ \hat{f}\left( u_{1}\right) -\hat{f}\left( u_{2}\right) %
\right] =\hat{f}^{\left( 3\right) }\left( u_{2}+\eta _{3}\left(
u_{1}-u_{2}\right) \right) \left( u_{1}-u_{2}\right) \left( D_{\xi
}u_{1}\right) ^{2}+\text{ }
\]%
\[
\hat{f}^{\left( 2\right) }\left( u_{2}\right) \left( D_{\xi }u_{1}-D_{\xi
}u_{2}\right) \left( D_{\xi }u_{1}+D_{\xi }u_{2}\right) +
\]%
\[
\hat{f}^{\left( 2\right) }\left( u_{2}+\eta _{4}\left( u_{1}-u_{2}\right)
\right) \left( u_{1}-u_{2}\right) D_{\xi }^{2}u_{1}+\hat{f}^{\left( 1\right)
}\left( u_{2}\right) \left( D_{\xi }^{2}u_{1}-D_{\xi }^{2}u_{2}\right) ,
\]%
where $0<\eta _{i}<1,$ $i=1,2,3,4.$ Thus, using H\"{o}lder's and Nirenberg's
inequality, we have%
\begin{equation}
\left\Vert \hat{f}\left( u_{1}\right) -\hat{f}\left( u_{2}\right)
\right\Vert _{X_{\infty }}\leq \bar{f}\left( M+1\right) \left\Vert
u_{1}-u_{2}\right\Vert _{X_{\infty }},  \tag{3.10}
\end{equation}%
\begin{equation}
\left\Vert \hat{f}\left( u_{1}\right) -\hat{f}\left( u_{2}\right)
\right\Vert _{X_{p}}\leq \bar{f}\left( M+1\right) \left\Vert
u_{1}-u_{2}\right\Vert _{X_{p}},  \tag{3.11}
\end{equation}%
\begin{equation}
\left\Vert D_{\xi }\left[ \hat{f}\left( u_{1}\right) -\hat{f}\left(
u_{2}\right) \right] \right\Vert _{X_{p}}\leq \left( M+1\right) \bar{f}%
\left( M+1\right) \left\Vert u_{1}-u_{2}\right\Vert _{X_{\infty }}+ 
\tag{3.12}
\end{equation}%
\[
\bar{f}\left( M+1\right) \left\Vert \hat{f}\left( u_{1}\right) -\hat{f}%
\left( u_{2}\right) \right\Vert _{X_{p}},
\]%
\[
\left\Vert D_{\xi }^{2}\left[ \hat{f}\left( u_{1}\right) -\hat{f}\left(
u_{2}\right) \right] \right\Vert _{X_{p}}\leq \left( M+1\right) \bar{f}%
\left( M+1\right) \left\Vert u_{1}-u_{2}\right\Vert _{X_{\infty }}\left\Vert
D_{\xi }^{2}u_{1}\right\Vert _{Y^{2,p}}^{2}+
\]%
\[
\bar{f}\left( M+1\right) \left\Vert D_{\xi }\left( u_{1}-u_{2}\right)
\right\Vert _{Y^{2,p}}\left\Vert D_{\xi }\left( u_{1}+u_{2}\right)
\right\Vert _{Y^{2,p}}+
\]%
\[
\bar{f}\left( M+1\right) \left\Vert u_{1}-u_{2}\right\Vert _{X_{\infty
}}\left\Vert D_{\xi }^{2}u_{1}\right\Vert _{X_{p}}+\bar{f}\left( M+1\right)
\left\Vert D_{\xi }\left( u_{1}-u_{2}\right) \right\Vert _{X_{p}}\leq 
\]%
\begin{equation}
C^{2}\bar{f}\left( M+1\right) \left\Vert u_{1}-u_{2}\right\Vert _{X_{\infty
}}\left\Vert u_{1}\right\Vert _{X_{\infty }}\left\Vert D_{\xi
}^{2}u_{1}\right\Vert _{X_{p}}+  \tag{3.13}
\end{equation}%
\[
C^{2}\bar{f}\left( M+1\right) \left\Vert u_{1}-u_{2}\right\Vert _{X_{\infty
}}^{\frac{1}{2}}\left\Vert D_{\xi }^{2}\left( u_{1}-u_{2}\right) \right\Vert
_{X_{p}}\left\Vert u_{1}+u_{2}\right\Vert _{X_{\infty }}^{\frac{1}{2}%
}\left\Vert D_{\xi }^{2}\left( u_{1}+u_{2}\right) \right\Vert _{X_{p}}
\]%
\[
+\left( M+1\right) \bar{f}\left( M+1\right) \left\Vert
u_{1}-u_{2}\right\Vert _{X_{\infty }}+\bar{f}\left( M+1\right) \left\Vert
D_{\xi }^{2}\left( u_{1}-u_{2}\right) \right\Vert _{X_{p}}\leq 
\]%
\[
3C^{2}\left( M+1\right) ^{2}\bar{f}\left( M+1\right) \left\Vert
u_{1}-u_{2}\right\Vert _{X_{\infty }}+2C^{2}\left( M+1\right) \bar{f}\left(
M+1\right) \left\Vert D_{\xi }^{2}\left( u_{1}-u_{2}\right) \right\Vert
_{X_{p}},
\]%
where $C$ is the constant in Lemma $3.1$. From $\left( 3.10\right) -\left(
3.11\right) $, using Minkowski's inequality for integrals, Fourier
multiplier theorems for operator-valued functions in $X_{p}$ spaces and
Young's inequality, we obtain%
\[
\left\Vert G\left( u_{1}\right) -G\left( u_{2}\right) \right\Vert _{Y\left(
T\right) }\leq \dint\limits_{0}^{t}\left\Vert u_{1}-u_{2}\right\Vert
_{X_{\infty }}d\tau +\dint\limits_{0}^{t}\left\Vert u_{1}-u_{2}\right\Vert
_{Y^{2,p}}d\tau +
\]%
\[
\dint\limits_{0}^{t}\left\Vert f\left( u_{1}\right) -f\left( u_{2}\right)
\right\Vert _{X_{\infty }}d\tau +\dint\limits_{0}^{t}\left\Vert f\left(
u_{1}\right) -f\left( u_{2}\right) \right\Vert _{Y^{2,p}}d\tau \leq 
\]%
\[
T\left[ 1+C_{1}\left( M+1\right) ^{2}\bar{f}\left( M+1\right) \right]
\left\Vert u_{1}-u_{2}\right\Vert _{Y\left( T\right) },
\]%
where $C_{1}$ is a constant. If $T$ satisfies $\left( 3.9\right) $ and the
following inequality holds 
\begin{equation}
T\leq \frac{1}{2}\left[ 1+C_{1}\left( M+1\right) ^{2}\bar{f}\left(
M+1\right) \right] ^{-1},  \tag{3.14}
\end{equation}%
then 
\[
\left\Vert Gu_{1}-Gu_{2}\right\Vert _{Y\left( T\right) }\leq \frac{1}{2}%
\left\Vert u_{1}-u_{2}\right\Vert _{Y\left( T\right) }.
\]

That is, $G$ is a constructive map. By contraction mapping principle we know
that $G(u)$ has a fixed point $u(x,t)\in $ $Q\left( M;T\right) $ that is a
solution of $(1.1)-(1.2)$. From $\left( 2.5\right) $ we get that $u$ is a
solution of the following integral equation 
\[
u\left( t,x\right) =S_{1}\left( t\right) \varphi \left( x\right)
+S_{2}\left( t\right) \psi \left( x\right) + 
\]%
\[
+\dint\limits_{0}^{t}F^{-1}\left[ S\left( t-\tau ,\xi \right) L\left( \xi
\right) \hat{f}\left( u\right) \left( \xi ,\tau \right) \right] d\tau ,\text{
}t\in \left( 0,T\right) . 
\]

Let us show that this solution is a unique in $Y\left( T\right) $. Let $%
u_{1} $, $u_{2}\in Y\left( T\right) $ are two solution of the problem $%
(1.1)-(1.2)$. Then%
\begin{equation}
u_{1}-u_{2}=\dint\limits_{0}^{t}F^{-1}\left[ S\left( t-\tau ,\xi \right)
L\left( \xi \right) \hat{f}\left( u_{1}\right) \left( \xi ,\tau \right) -%
\hat{f}\left( u_{2}\right) \left( \xi ,\tau \right) \right] d\tau . 
\tag{3.15}
\end{equation}
By the definition of the space $Y\left( T\right) $, we can assume that%
\[
\left\Vert u_{1}\right\Vert _{X_{\infty }}\leq C_{1}\left( T\right) ,\text{ }%
\left\Vert u_{1}\right\Vert _{X_{\infty }}\leq C_{1}\left( T\right) . 
\]

Hence, by Minkowski's inequality for integrals and Theorem 2.2 we obtain
from $\left( 3.15\right) $

\begin{equation}
\left\Vert u_{1}-u_{2}\right\Vert _{Y^{2,p}}\leq C_{2}\left( T\right) \text{ 
}\dint\limits_{0}^{t}\left\Vert u_{1}-u_{2}\right\Vert _{Y^{2,p}}d\tau . 
\tag{3.16}
\end{equation}

From $(3.16)$ and Gronwall's inequality, we have $\left\Vert
u_{1}-u_{2}\right\Vert _{Y^{2,p}}=0$, i.e. problem $(1.1)-(1.2)$ has a
unique solution which belongs to $Y\left( T\right) .$ That is, we obtain the
first part of the assertion. Now, let $\left[ 0\right. ,\left. T_{0}\right) $
be the maximal time interval of existence for $u\in Y\left( T_{0}\right) $.
It remains only to show that if $(3.3)$ is satisfied, then $T_{0}=\infty $.
Assume contrary that, $\left( 3.3\right) $ holds and $T_{0}<\infty .$ For $%
T\in \left[ 0\right. ,\left. T_{0}\right) ,$ we consider the following
integral equation

\begin{equation}
\upsilon \left( x,t\right) =S_{1}\left( t\right) u\left( x,T\right)
+S_{2}\left( t\right) u_{t}\left( x,T\right) +  \tag{3.17}
\end{equation}

\[
\dint\limits_{0}^{t}F^{-1}\left[ S\left( t-\tau ,\xi \right) L\left( \xi
\right) \hat{f}\left( \upsilon \right) \left( \xi ,\tau \right) \right]
d\tau ,\text{ }t\in \left( 0,T\right) . 
\]%
By virtue of $(3.3)$, for $T^{\prime }>T$ we have 
\[
\sup_{t\in \left[ 0\right. ,\left. T\right) }\left( \left\Vert u\right\Vert
_{Y^{2,p}}+\left\Vert u\right\Vert _{X_{\infty }}+\left\Vert
u_{t}\right\Vert _{Y^{2,p}}+\left\Vert u_{t}\right\Vert _{X_{\infty
}}\right) <\infty . 
\]

By reasoning as a first part of theorem and by contraction mapping
principle, there is a $T^{\ast }\in \left( 0,T_{0}\right) $ such that for
each $T\in \left[ 0\right. ,\left. T_{0}\right) ,$ the equation $\left(
3.17\right) $ has a unique solution $\upsilon \in Y\left( T^{\ast }\right) .$
The estimates $\left( 3.9\right) $ and $\left( 3.14\right) $ imply that $%
T^{\ast }$ can be selected independently of $T\in \left[ 0\right. ,\left.
T_{0}\right) .$ Set $T=T_{0}-\frac{T^{\ast }}{2}$ and define 
\begin{equation}
\tilde{u}\left( x,t\right) =\left\{ 
\begin{array}{c}
u\left( x,t\right) ,\text{ }t\in \left[ 0,T\right]  \\ 
\upsilon \left( x,t-T\right) \text{, }t\in \left[ T,T_{0}+\frac{T^{\ast }}{2}%
\right] 
\end{array}%
\right. .  \tag{3.18}
\end{equation}

By construction $\tilde{u}\left( x,t\right) $ is a solution of the problem $%
(1.1)-(1.2)$ on $\left[ T,T_{0}+\frac{T^{\ast }}{2}\right] $ and in view of
local uniqueness, $\tilde{u}\left( x,t\right) $ extends $u.$ This is against
to the maximality of $\left[ 0\right. ,\left. T_{0}\right) $, i.e we obtain $%
T_{0}=\infty .$

\begin{center}
\textbf{4. Applications}\qquad
\end{center}

\bigskip In this section we give some application of Theorem 3.1.

\textbf{1. }Let%
\[
L_{0}=L_{1}=L_{2}=A_{1}=\dsum\limits_{\left\vert \alpha \right\vert \leq
2}a_{\alpha }D^{\alpha }, 
\]%
where $a_{\alpha }$ are complex numbers.

Then the problem $\left( 1.1\right) -\left( 1.2\right) $ is reduced to the
Cauchy problem for the following Boussinesq equation%
\begin{equation}
u_{tt}+A_{1}u_{tt}+A_{1}u=A_{1}f\left( x,t,u\right) ,\text{ }x\in R^{2},%
\text{ }t\in \left( 0,T\right) ,  \tag{4.1}
\end{equation}%
\[
u\left( x,0\right) =\varphi \left( x\right) +\dint\limits_{0}^{T}\alpha
\left( \sigma \right) u\left( x,\sigma \right) d\sigma ,\text{ }u_{t}\left(
x,0\right) =\psi \left( x\right) +\dint\limits_{0}^{T}\beta \left( \sigma
\right) u_{t}\left( x,\sigma \right) d\sigma ,
\]%
here 
\[
\varphi ,\text{ }\psi \in W^{s,p}\left( R^{2}\right) \text{, }s>\frac{2}{p},%
\text{ }p\in \left( 1,\infty \right) .
\]

Assume 
\[
\text{ }A_{1}\left( \xi \right) =\dsum\limits_{\left\vert \alpha \right\vert
\leq 2}a_{\alpha }\xi _{1}^{\alpha _{1}}\xi _{2}^{\alpha _{2}}>0\text{ for }%
\xi =\left( \xi _{1},\xi _{2}\right) \in R^{2}
\]%
\ Then it is not hard to see that%
\[
A_{1}\left( \xi \right) \neq 0,-1\text{ and}\left\vert Q^{-\frac{1}{2}%
}\left( \xi \right) \right\vert \leq 1\text{ for }\xi \in R^{2}\text{,}
\]%
where%
\[
Q\left( \xi \right) =A_{1}\left( \xi \right) \left[ 1+A_{1}\left( \xi
\right) \right] ^{-1}.
\]%
Hence, the Condition 2.1 is satisfied. Let 
\[
X_{p}=L^{p}\left( R^{2}\right) \text{, }1\leq p\leq \infty ,\text{ }%
Y^{s,p}=L^{s,p}\left( R^{2}\right) .
\]

Hence, from Theorem 3.1 we obtain:

\textbf{Theorem 4.1. }Assume that the function $u\rightarrow $ $f\left(
x,t,u\right) $: $R^{2}\times \left[ 0,T\right] \times B_{p}^{2-\frac{1}{p}%
}\left( R^{2}\right) \rightarrow L^{p}\left( R^{2}\right) $ is measurable in 
$\left( x,t\right) \in R^{2}\times \left[ 0,T\right] $ for $u\in B_{p}^{2-%
\frac{1}{p}}\left( R^{2}\right) .$ Moreover, $f\left( x,t,u\right) $ is
continuous in $u\in B_{p}^{2-\frac{1}{p}}\left( R^{2}\right) $ and%
\[
f\left( x,t,u\right) \in C^{\left( 3\right) }\left( B_{p}^{2-\frac{1}{p}%
}\left( R^{2}\right) ;\mathbb{C}\right) 
\]%
uniformly with respect to $\left( x,t\right) \in R^{2}\times \left[ 0,T%
\right] $. Then for $\varphi ,$ $\psi $ $\in Y_{\infty }^{2,p}$ and $p\in
\left( 1,\infty \right) $ problem $\left( 4.1\right) $ has a unique local
strange solution $u\in C^{\left( 2\right) }\left( \left[ 0,\right. \left.
T_{0}\right) ;Y_{\infty }^{2,p}\right) $, where $T_{0}$ is a maximal time
interval that is appropriately small relative to $M$. Moreover, if

\[
\sup_{t\in \left[ 0\right. ,\left. T_{0}\right) }\left( \left\Vert
u\right\Vert _{Y^{2,p}}+\left\Vert u\right\Vert _{X_{\infty }}+\left\Vert
u_{t}\right\Vert _{Y^{2,p}}+\left\Vert u_{t}\right\Vert _{X_{\infty
}}\right) <\infty 
\]%
then $T_{0}=\infty .$

\textbf{2. }Let%
\[
L_{0}=L_{1}=L_{2}=A_{2}=\dsum\limits_{\left\vert \alpha \right\vert \leq
4}a_{\alpha }D^{\alpha }, 
\]%
where $a_{\alpha }$ are complex numbers, $\alpha =\left( \alpha _{1},\alpha
_{2},\alpha _{3}\right) ,$ $\alpha _{k}$ are natural numbers and $\left\vert
\alpha \right\vert =$ $\dsum\limits_{k=1}^{3}$ $\alpha _{k}.$

Then the problem $\left( 1.1\right) -\left( 1.2\right) $ is reduced to the
Cauchy problem for the following Boussinesq equation%
\begin{equation}
u_{tt}+A_{2}u_{tt}+A_{2}u=A_{2}f\left( x,t,u\right) ,\text{ }x\in R^{3},%
\text{ }t\in \left( 0,T\right) ,  \tag{4.2}
\end{equation}%
\[
u\left( x,0\right) =\varphi \left( x\right) +\dint\limits_{0}^{T}\alpha
\left( \sigma \right) u\left( x,\sigma \right) d\sigma ,\text{ }u_{t}\left(
x,0\right) =\psi \left( x\right) +\dint\limits_{0}^{T}\beta \left( \sigma
\right) u_{t}\left( x,\sigma \right) d\sigma ,
\]%
where 
\[
\varphi ,\text{ }\psi \in W^{s,p}\left( R^{3}\right) \text{, }s>\frac{3}{p},%
\text{ }p\in \left( 1,\infty \right) .
\]%
Assume 
\[
A_{2}\left( \xi \right) =\dsum\limits_{\left\vert \alpha \right\vert \leq
4}a_{\alpha }\xi _{1}^{\alpha _{1}}\xi _{2}^{\alpha _{2}}\xi _{3}^{\alpha
_{3}}\neq 0,\text{ }-1\text{ for all }\xi =\left( \xi _{1},\xi _{2},\xi
_{3}\right) \in R^{3}.
\]%
\ Then it is not hard to see that, there exists a positive constant $M$ such
that%
\[
\left\vert Q^{\frac{1}{2}}\left( \xi \right) \right\vert \leq 1\text{ for }%
\xi \in R^{3}\text{ and }p\in \left( 1,\infty \right) ,
\]%
where 
\[
Q\left( \xi \right) =A_{2}\left( \xi \right) \left[ 1+A_{2}\left( \xi
\right) \right] ^{-1}.
\]%
Therefore, the Condition 2.1 is satisfied.

Let 
\[
X_{p}=L^{p}\left( R^{3}\right) \text{, }1\leq p\leq \infty ,\text{ }%
Y^{s,p}=L^{s,p}\left( R^{3}\right) . 
\]

Hence, from Theorem 3.1 we obtain:

\textbf{Theorem 4.2. }Suppose that the function $u\rightarrow $ $f\left(
x,t,u\right) $: $R^{3}\times \left[ 0,T\right] \times B_{p}^{2-\frac{1}{p}%
}\left( R^{3}\right) \rightarrow L^{p}\left( R^{3}\right) $ is measurable in 
$\left( x,t\right) \in R^{3}\times \left[ 0,T\right] $ for $u\in B_{p}^{2-%
\frac{1}{p}}\left( R^{3}\right) .$ Moreover, $f\left( x,t,u\right) $ is
continuous in $u\in B_{p}^{2-\frac{1}{p}}\left( R^{3}\right) $ and%
\[
f\left( x,t,u\right) \in C^{\left( 3\right) }\left( B_{p}^{2-\frac{1}{p}%
}\left( R^{3}\right) ;\mathbb{C}\right) 
\]%
uniformly with respect to $\left( x,t\right) \in R^{3}\times \left[ 0,T%
\right] $. Then for $\varphi ,$ $\psi $ $\in Y_{\infty }^{2,p}$ and $p\in
\left( 1,\infty \right) $ problem $\left( 4.1\right) $ has a unique local
strange solution%
\[
u\in C^{\left( 2\right) }\left( \left[ 0,\right. \left. T_{0}\right)
;Y_{\infty }^{2,p}\right) , 
\]%
where $T_{0}$ is a maximal time interval that is appropriately small
relative to $M$. Moreover, if

\[
\sup_{t\in \left[ 0\right. ,\left. T_{0}\right) }\left( \left\Vert
u\right\Vert _{Y^{2,p}}+\left\Vert u\right\Vert _{X_{\infty }}+\left\Vert
u_{t}\right\Vert _{Y^{2,p}}+\left\Vert u_{t}\right\Vert _{X_{\infty
}}\right) <\infty 
\]%
then $T_{0}=\infty .$

\textbf{3. }Let%
\[
L_{0}=\dsum\limits_{\left\vert \alpha \right\vert \leq 4}a_{0\alpha
}D^{\alpha }\text{, }L_{1}=\dsum\limits_{\left\vert \alpha \right\vert \leq
2}a_{1\alpha }D^{\alpha },\text{ }L_{2}=\dsum\limits_{\left\vert \alpha
\right\vert \leq 4}a_{2\alpha }D^{\alpha }, 
\]%
where $a_{\alpha i}$ are complex numbers, $\alpha =\left( \alpha _{1},\alpha
_{2},\alpha _{3}\right) ,$ $\alpha _{k}$ are natural numbers and $\left\vert
\alpha \right\vert =$ $\dsum\limits_{k=1}^{3}$ $\alpha _{k}.$

Then the problem $\left( 1.1\right) -\left( 1.2\right) $ is reduced to
Cauchy problem for the following Boussinesq equation%
\begin{equation}
u_{tt}+L_{0}u_{tt}+L_{1}u=L_{2}f\left( x,t,u\right) ,\text{ }x\in R^{3},%
\text{ }t\in \left( 0,T\right) ,  \tag{4.3}
\end{equation}%
\[
u\left( x,0\right) =\varphi \left( x\right) +\dint\limits_{0}^{T}\alpha
\left( \sigma \right) u\left( x,\sigma \right) d\sigma ,\text{ }u_{t}\left(
x,0\right) =\psi \left( x\right) +\dint\limits_{0}^{T}\beta \left( \sigma
\right) u_{t}\left( x,\sigma \right) d\sigma ,
\]%
where 
\[
\varphi ,\text{ }\psi \in W^{s,p}\left( R^{3}\right) \text{, }s>\frac{3}{p},%
\text{ }p\in \left( 1,\infty \right) .
\]

\bigskip\ Assume 
\[
\text{ }L_{0}\left( \xi \right) >0\text{, }L_{1}\left( \xi \right) >0\text{
for }\xi =\left( \xi _{1},\xi _{2},\xi _{3}\right) \in R^{3}. 
\]

Since $m_{0}-m_{1}=m_{2}-m_{1}=2,$ for $s\geq 1+\frac{3}{p}$ the Condition
2.1 is satisfied.

Let 
\[
X_{p}=L^{p}\left( R^{3}\right) \text{, }1\leq p\leq \infty ,\text{ }%
Y^{s,p}=L^{s,p}\left( R^{3}\right) . 
\]

Hence, from Theorem 3.1 we obtain:

\textbf{Theorem 4.3. }Suppose that the function $u\rightarrow $ $f\left(
x,t,u\right) $: $R^{3}\times \left[ 0,T\right] \times B_{p}^{2-\frac{1}{p}%
}\left( R^{3}\right) \rightarrow L^{p}\left( R^{3}\right) $ is measurable in 
$\left( x,t\right) \in R^{3}\times \left[ 0,T\right] $ for $u\in B_{p}^{2-%
\frac{1}{p}}\left( R^{3}\right) .$ Moreover, $f\left( x,t,u\right) $ is
continuous in $u\in B_{p}^{2-\frac{1}{p}}\left( R^{3}\right) $ and%
\[
f\left( x,t,u\right) \in C^{\left( 3\right) }\left( B_{p}^{2-\frac{1}{p}%
}\left( R^{3}\right) ;\mathbb{C}\right) 
\]%
uniformly with respect to $\left( x,t\right) \in R^{3}\times \left[ 0,T%
\right] $. Then for $\varphi ,$ $\psi $ $\in Y_{\infty }^{2,p}$ and $s\geq 1+%
\frac{3}{p},$ $p\in \left( 1,\infty \right) $ problem $\left( 4.3\right) $
has a unique local strange solution%
\[
u\in C^{\left( 2\right) }\left( \left[ 0,\right. \left. T_{0}\right)
;Y_{\infty }^{2,p}\right) , 
\]%
where $T_{0}$ is a maximal time interval that is appropriately small
relative to $M$. Moreover, if

\[
\sup_{t\in \left[ 0\right. ,\left. T_{0}\right) }\left( \left\Vert
u\right\Vert _{Y^{2,p}}+\left\Vert u\right\Vert _{X_{\infty }}+\left\Vert
u_{t}\right\Vert _{Y^{2,p}}+\left\Vert u_{t}\right\Vert _{X_{\infty
}}\right) <\infty 
\]%
then $T_{0}=\infty .$

\textbf{References}

\begin{quote}
\ \ \ \ \ \ \ \ \ \ \ \ \ \ \ \ \ \ \ \ \ \ \ \ 
\end{quote}

\begin{enumerate}
\item V.G. Makhankov, Dynamics of classical solutions (in non-integrable
systems), Phys. Lett. C 35, (1978) 1--128.

\item G.B. Whitham, Linear and Nonlinear Waves, Wiley--Interscience, New
York, 1975.

\item N.J. Zabusky, Nonlinear Partial Differential Equations, Academic
Press, New York, 1967.

\item C. G. Gal and A. Miranville, Uniform global attractors for
non-isothermal viscous and non-viscous Cahn--Hilliard equations with dynamic
boundary conditions, Nonlinear Analysis: Real World Applications 10 (2009)
1738--1766.

\item T. Kato, T. Nishida, A mathematical justification for Korteweg--de
Vries equation and Boussinesq equation of water surface waves, Osaka J.
Math. 23 (1986) 389--413.

\item A. Clarkson, R.J. LeVeque, R. Saxton, Solitary-wave interactions in
elastic rods, Stud. Appl. Math. 75 (1986) 95--122.

\item P. Rosenau, Dynamics of nonlinear mass-spring chains near continuum
limit, Phys. Lett. 118A (1986) 222--227.

\item S. Wang, G. Chen, Small amplitude solutions of the generalized IMBq
equation, J. Math. Anal. Appl. 274 (2002) 846--866.

\item S. Wang, G. Chen,The Cauchy Problem for the Generalized IMBq equation
in $W^{s,p}\left( R^{n}\right) $, J. Math. Anal. Appl. 266, 38--54 (2002).

\item J.L. Bona, R.L. Sachs, Global existence of smooth solutions and
stability of solitary waves for a generalized Boussinesq equation, Comm.
Math. Phys. 118 (1988) 15--29.

\item F. Linares, Global existence of small solutions for a generalized
Boussinesq equation, J. Differential Equations 106 (1993) 257--293.

\item Y. Liu, Instability and blow-up of solutions to a generalized
Boussinesq equation, SIAM J. Math. Anal. 26 (1995) 1527--1546.

\item S. Piskarev and S.-Y. Shaw, Multiplicative perturbations of semigroups
and applications to step responses and cumulative outputs, J. Funct. Anal.
128 (1995), 315-340.

\item N. Kutev, N. Kolkovska, and M. Dimova, \textquotedblleft Global
existence of Cauchy problem for Boussinesq paradigm
equation,\textquotedblright\ Computers and Mathematics with Applications,
65(3) (2013), 500--511,

\item S. Lai, Y.H. Wu, The asymptotic solution of the Cauchy problem for a
generalized Boussinesq equation, Discrete Contin. Dyn. Syst. Ser. B 3 (2003).

\item Girardi, M., Lutz, W., Operator-valued Fourier multiplier theorems on $%
L_{p}(X)$ and geometry of Banach spaces, J. Funct. Anal., 204(2), 320--354,
2003.

\item H. Triebel, Interpolation theory, Function spaces, Differential
operators, North-Holland, Amsterdam, 1978.

\item H. Triebel, Fractals and spectra, Birkhauser Verlag, Related to
Fourier analysis and function spaces, Basel, 1997.

\item L. Nirenberg, On elliptic partial differential equations, Ann. Scuola
Norm. Sup. Pisa 13 (1959), 115--162.

\item S. Klainerman, Global existence for nonlinear wave equations, Comm.
Pure Appl. Math. 33 (1980), 43--101.

\item A. Ashyralyev, N. Aggez, Nonlocal boundary value hyperbolic problems
involving Integral conditions, Bound.Value Probl., 2014 V. 2014:214.

\item L. S. Pulkina, A non local problem with integral conditions for
hyperbolice quations, Electron.J.Differ.Equ.1999,45 (1999)
\end{enumerate}

\end{document}